\documentclass[preprint,12pt]{elsarticle}

\usepackage{amssymb}
\usepackage{amsmath}

\newtheorem{thm}{Theorem}

\journal{}

\begin{document}

\begin{frontmatter}



\title{A note on ‘‘The asymptotic uniform distribution of subset 
sums’’}


\author{Yilong Hu\footnote{ email address: huyl10@sjtu.edu.cn}} 

\affiliation{organization={Shanghai Jiao Tong University},
            addressline={800 Dongchuan Rd.}, 
            city={Shanghai},
            postcode={200240},
            country={China}}

\begin{abstract}
We find out that the main result of the article ‘‘The asymptotic uniform distribution of subset sums’’ can be proven much more easily, using an explicit formula proposed by Li and Wan \cite{li2012counting}.
\end{abstract}


\end{frontmatter}



Let $G$ be a finite abelian group of order $n$. It's well known that there exists positive integers $n_1,\dots,n_s$ such that $G \cong \mathbb{Z}_{n_1}\times \mathbb{Z}_{n_2}\times \dots \times \mathbb{Z}_{n_s}$.  For any $b\in G$ (which can be recognized as $(b_1,\dots,b_s)$), Let $N(k, b)$ be the number of $k$-subsets of $A$ whose elements sum to $b$. In \cite{wang2026asymptotic}, author proved the following theorem:
\begin{thm}
     Let $G$ be a finite abelian group of order $n$ and let $k := k(n)$ be a function of $n$. Then
$$ \lim_{n \to \infty} \frac{\min_{b \in G} |N(k, b)|}{\max_{b \in G} |N(k, b)|} = 1 $$
holds for all $4 \le k \le \lfloor \frac{n}{2} \rfloor + 1$.
\end{thm}

The goal of this note is to give a shorter proof of this theorem. Firstly, the following explicit formula of $N(k, b)$ is introduced:
\begin{thm} 
$$ N(k, b) = \frac{1}{n} \sum_{r | (n, k)} (-1)^{k + \frac{k}{r}} \binom{n/r}{k/r} \Phi(r, b), $$
    where $\Phi(r, b) = \sum_{d|r, (n_i, d)|b_i} \mu(r/d) \prod_{i=1}^s (n_i, d)$ and $\mu$ is the usual M\"{o}bius function defined over the integers.
\end{thm}

Note that the $r=1$ term is exactly $\binom{n}{k}/n$, while the others terms have a trivial upper bound $k\binom{n/2}{k/2}$ (here one may assume that both $n$ and $k$ are even.) One only need to prove that 
$$\lim_{n \to \infty}\frac{nk^2\binom{n/2}{k/2}}{\binom{n}{k}}=0. $$

 It's easy to show that 

$${n/2 \choose k/2}/{n\choose k} \le (\frac{k}{n})^{k/2}.$$

Let $L_n(k) = \log (nk^2(\frac{k}{n})^{k/2})$, Then $\frac{D^2L_n}{Dk^2}=\frac{1}{2k}-\frac{2}{k^2}$, which is always positive when $k\ge 4$. Consequentially, the maximal value is taken only at $k=4$ or $k=\frac{n}{2}+1$.

When $k=4$, $L_n(k) = \log (16n(\frac{4}{n})^{2})=\log (256/n)$; when $k=\frac{n}{2}+1$, $L_n(k) = \log (n(\frac{n+2}{2})^2(\frac{n+2}{2n})^{\frac{n+2}{4}})$. It is evident that in both cases, $L_n(k)$ approaches negative infinity as $n$ tends to infinity. This finishes the proof.
\bibliographystyle{elsarticle-num} 
\bibliography{name}






\end{document}